\documentclass[a4paper,12pt]{article}
\RequirePackage{amsmath,amsthm,amssymb}
\RequirePackage{hyperref}
\RequirePackage{graphicx}

\newcommand{\half}{\textstyle{\frac{1}{2}}}
\newcommand{\dd}{\mathrm{d}}

\begin{document}

\begin{Large}
\centerline{\bf Derangements and Continued Fractions for $e$}
\end{Large}
\medskip
\centerline{\sc Peter Lynch\footnote{Peter.Lynch@ucd.ie}}
\centerline{School of Mathematics \&\ Statistics, University College Dublin}

\begin{abstract}
Several continued fraction expansions for $e$ have been
produced by an automated conjecture generator (ACG) called
\emph{The Ramanujan Machine}. Some of these were already
known, some have recently been proved and some remain
unproven. While an ACG can produce interesting putative
results, it gives very limited insight into their
significance. In this paper, we derive an elegant continued
fraction expansion, equivalent to a result from 
the Ramanujan Machine, using the sequence of ratios of
factorials to subfactorials or derangement numbers.
\end{abstract}



\subsection*{Arrangements and Derangements}

Six students entering an examination hall place their cell-phones in a box.
After the exam, they each grab a phone at random as they rush out.
What is the likelihood that none of them gets their own phone?
The surprising answer --- about 37\% whatever the number of students ---
emerges from the theory of derangements.


We may call any permutation of the elements of a set an arrangement.
A \emph{derangement} is an arrangement for which every element is moved
from its original position.
Thus, a derangement is a permutation that has no \emph{fixed points.}
The number of derangements of a set of $n$ elements is also called the
\emph{subfactorial} of $n$. Various notations are used for subfactorials:
$!n$, $d_n$ and $n$\textexclamdown\ are common; we will use $!n$
(read as `bang-en').

Derangements were first considered by Pierre Reymond de Montmort.
In 1713, with help from Nicholas Bernoulli, he managed to find an
expression for the connection between $!n$ and $n!$.  The answer,
which he obtained using the inclusion-exclusion principle
(Zeilberger, 2008, pg.~560), is
\begin{equation}
!n = n! \left(1 - \frac{1}{1!} + \frac{1}{2!} - \frac{1}{3!} 
                + \frac{1}{4!} - \cdots \pm\frac{1}{n!}\right)
= n! \sum_{k=0}^n \frac{(-1)^k}{k!} \,.
\label{eq:deMont}
\end{equation}
Of course, we see from this that $\lim_{n\to\infty} (!n) = n!/e$.
In fact, we can write a more precise connection
between derangements and arrangements:
$$
!n = \left\lfloor\frac{n!+\half}{e}\right\rfloor \,.
$$
This means that $!n$ is the nearest whole number to $n!/e$. 


The number $!n$ of derangements of an $n$-element set 
may be calculated using a second-order recurrence relation:
$$
 !n = (n-1) [!(n-1)+!(n-2)]
$$
with $!0 = 1$ and $!1 = 0$.
The subfactorials also satisfy a first-order recurrence relation,
$$
!n = n\times !(n-1) + (-1)^n \,, \qquad [\mbox{compare\ \ } n! = n\times(n-1)! ]
$$
with initial condition $!0 = 1$.  The first eight values of $!n$ are
$\{ 1, 0, 1, 2, 9, 44, 265, 1\,854  \}$. 

There is an integral expression for the subfactorial,
\begin{equation}
!n = \int _{0}^\infty (x-1)^n e^{-x}\dd x \,,
\qquad \left[
\mbox{compare\ \ } n! = \int _{0}^\infty x^n e^{-x}\dd x \right] \,.
\label{eq:bangInt}
\end{equation}
Expansion of (\ref{eq:bangInt}) yields de~Montmort's result (\ref{eq:deMont}).
It also allows extension of the subfactorial function to non-integral arguments
($!x$) and analytic continuation to the complex plane ($!z$). 


\subsection*{Continued Fractions and Convergents}

The continued fraction expansion of an irrational number $x$ is written,
in expanded form (centre) and concise form (right), as
$$
x = a_0 + \cfrac{1}{ a_1 + \cfrac{1}{ a_2 + \cfrac{1}{ a_3 + \ddots } }}
= [ a_0 ; a_1 , a_2 , a_3 , \dots ] 
$$
where $a_n$ are integers.
If $a_n$ is positive for $n \ge 1$ this is called a \emph{simple} continued fraction.

The \emph{generalized} continued fraction expansion is written
$$
x = b_0 + \cfrac{a_1}{ b_1 + \cfrac{a_2}{ b_2 + \cfrac{a_3}{ b_3 + \ddots } }} = 
b_0 +\frac{a_1}{b_1+}\,\frac{a_2}{b_2+}\,\frac{a_3}{b_3+}\,\frac{a_4}{b_4+}\,\cdots \,.
$$
We assume that $a_n$ and $b_n$ are integers. Truncating the expansion at various points,
we obtain the \emph{convergents}
$$
r_n = \frac{p_n}{q_n} =
b_0 +\frac{a_1}{b_1+}\,\frac{a_2}{b_2+}\,\frac{a_3}{b_3+}\,\frac{a_4}{b_4+}\,\cdots\,\frac{a_n}{b_n}
$$
where the numerators and denominators, $p_n$ and $q_n$, are integers.  We define the starting values
$$
p_{-1} = 1\,, \qquad  q_{-1} = 0\,, \qquad  p_0 = b_0\,, \qquad q_0 = 1 \,.
$$
Then, $p_{k}$ and $q_{k}$ for $k\ge 1$ are given by recurrence relations:
\begin{equation}
p_{k} = b_{k} p_{k-1} + a_{k} p_{k-2} \,, \qquad
q_{k} = b_{k} q_{k-1} + a_{k} q_{k-2} \,,
\label{eq:convergents}
\end{equation}
which may be proved by induction (Jones \&\ Thron, 1980, Pg.~20). 

This process can be inverted:
given a sequence of numerators $p_n$ and denominators $q_n$
(or just their ratios, the convergents $r_n = p_n/q_n$),
we can solve (\ref{eq:convergents}) for $a_n$ and $b_n$:
\begin{equation}
a_n = \frac{p_{n-1}q_{n} - p_{n}q_{n-1}} {p_{n-1}q_{n-2} - p_{n-2}q_{n-1}} \,, \qquad
b_n = \frac{p_{n}q_{n-2} - p_{n-2}q_{n}} {p_{n-1}q_{n-2} - p_{n-2}q_{n-1}}
\label{eq:numdenom}
\end{equation}
together with the starting values $b_0=p_0$, $a_1 = (p_1-b_0q_1)$ and $b_1=q_1$.

%

\subsection*{Continued Fractions for $e$}

Euler's number is usually defined as the limit
$ e = \lim_{n\to\infty}(1+1/n)^n$, which is the limit of the sequence
$$
\left\{ \frac{2^1}{1^1}, \frac{3^2}{2^2}, \frac{4^3}{3^3}, \dots ,\frac{(n+1)^n}{n^n}, \dots  \right\}
$$
The terms may be regarded as the convergents of a continued fraction, 
$$
r_n = \frac{p_n}{q_n} \,, \qquad\mbox{where}\qquad 
p_n = (n+1)^n \quad\mbox{and}\quad q_n = n^n \,.
$$
We can generate a continued fraction by using (\ref{eq:numdenom}). It begins as
\begin{equation}
1 +\frac{1}{1-}\,\frac{1}{5-}\,\frac{13}{10-}\,\frac{491}{196-}\,
\frac{487903}{9952-}\,\frac{2384329879}{958144-}\, \cdots \,.
\label{eq:cf1}
\end{equation}
The error of this expansion ($\log_{10}[r_n-e]$) as a function of truncation
 is shown in Fig.~\ref{fig:cfconv} (dashed line).
It is clear that the convergence is very slow.

Euler made extensive studies of continued fractions.
For example, his 50-page paper, \emph{Observations on continued fractions} (Euler, 1750),
contains numerous original results. One of his best-known expansions is
\begin{equation}
e = [2; 1,2,1,1,4,1,1,6,1,1,8,\dots]
\label{eq:EulerCF}
\end{equation}
The error of Euler's expansion is shown in Fig.~\ref{fig:cfconv} (dotted line).
It converges much faster than (\ref{eq:cf1}). There is a clear signal of period 3,
consistent with the recurring pattern $(1, 1, n)$ in (\ref{eq:EulerCF}).

\subsection*{Continued fraction from derangement numbers}

\begin{figure}[h]
\begin{center}
\includegraphics[width=0.85\textwidth]{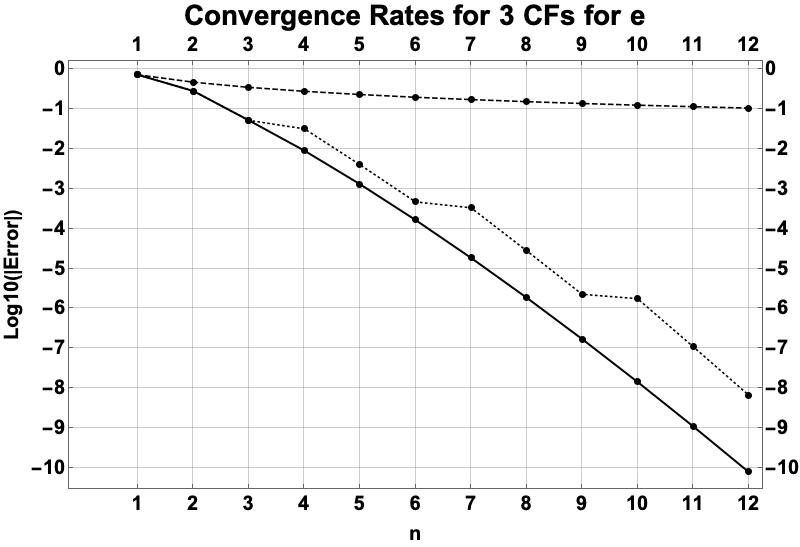}
\caption{Logarithm of the error $\log_{10}| r_n - e |$ in the continued
fraction expansions for $e$.
Dashed line: $r_n=(1+1/n)^n$, Eq.~(\ref{eq:cf1}) .
Dotted line: Convergents of Euler's expansion (\ref{eq:EulerCF}).
Solid line: $r_n=(n+1)!/!(n+1)$, Eq.~(\ref{eq:cf2}).}
\label{fig:cfconv}
\end{center}
\end{figure}

A beautiful continued fraction emerges from the relationship
between arrangements and derangements. We saw above that
$$
\left[
\frac{\mbox{Arrangements of $n$ elements}}
     {\mbox{Derangements of $n$ elements}}  \right]
= \frac{n!}{!n} \to e
$$
If we define the numerators and denominators of convergents to be
$$
p_n = n!  \quad\mbox{and}\quad q_n = \,!n \,,
$$
we can solve for the factors $a_n$ and $b_n$.
The starting values $p_0=1, p_1=1, q_0=1, q_1=0$ yield
$a_0=0, b_0=1, a_1=1, b_1=0$.  Then (\ref{eq:numdenom}) may be solved to yield
$a_n = b_n = (n-1)$ for $n\ge 2$. Thus we get the expansion
$$
e = 1 +\frac{1}{0+}\,\frac{1}{1+}\,\frac{2}{2+}\,\frac{3}{3+}\,\frac{4}{4+}\,\cdots \,.
$$
A small adjustment enables us to write this in the elegant form
\begin{equation}
e = 2 +\frac{2}{2+}\,\frac{3}{3+}\,\frac{4}{4+}\,\frac{5}{5+}\,\frac{6}{6+}\,\cdots \,.
\label{eq:cf2}
\end{equation}
The error of (\ref{eq:cf2}) is shown in Fig.~\ref{fig:cfconv} (solid line).
Convergence is more rapid than for the other two expansions.

\subsection*{The Ramanujan Machine}

An Automated Conjecture Generator (ACG) called \emph{The Ramanujan Machine}%
\footnote{G.~H.~Hardy, in his Introduction to Ramanujan's \emph{Collected Papers}
(1927), wrote that Ramanujan's mastery of of continued fractions was 
``beyond that of any mathematician in the world''.}
has been implemented by a team of mathematicians at the Israel Institute of
Technology.  This ACG system is capable of producing conjectures about mathematical
(and physical) constants, expressed in the form of continued fractions, using only 
numerical data as input.  A paper describing the system is available on the arXiv
preprint server (Raayoni, et al., 2020).

The Ramanujan Machine comprises algorithms designed to discover new conjectures,
running on a network of computers.
The goal of the project is to formulate conjectures that may then be proved mathematically.
The ACG has already generated a number of very interesting new conjectures,
as well as reproducing several results that were already well known.
The website (\url{http://www.ramanujanmachine.com/}) enables researchers to submit proofs
of conjectures, code new algorithms and (if they wish) allow access to their computers
for distributed computation.

While the Ramanujan Machine generates conjectures but not proofs, it has inspired a
complementary project using \emph{symbolic} rather than numerical computation.
Dougherty-Bliss and Zeilberger (2020) describe a system that generates automatic proofs 
of continued fraction expansions. Their system produced some infinite families of expansions
together with rigorous proofs of their validity.

One of the continued fractions discovered by the Ramanujan Machine is
\begin{equation}
\frac{1}{e-1} = \frac{1}{1+}\,\frac{2}{2+}\,\frac{3}{3+}\,\frac{4}{4+}\,\frac{5}{5+}\,\frac{6}{6+}\,\cdots \,,
\label{eq:OneOnEminus1}
\end{equation}
which is easily seen to be equivalent to (\ref{eq:cf2}) above. This is indicated in
Raayoni, et al.~(2020) as a ``known'' result. 
A proof was presented by Kadyrov and Mashurov (2019).  Lu (2019) gave
elementary proofs of other generalized continued fraction formulae for $e$.
However, the connection with derangement numbers was not made by any of these authors.

Dougherty-Bliss and Zeilberger (2020) proved a generalized expansion of which
the Ramanujan Machine result is a special case.
They noted the occurrence of derangement numbers in their expansion,
describing this as a ``remarkable coincidence'', and further commenting
that ``There does not seem to be any immediate combinatorial reason
for the derangement numbers to appear.''  Our above derivation of
(\ref{eq:cf2}), starting from the
ratio of factorials to subfactorials, makes the connection clear.




\end{document}